\documentclass{letter}

\input{xy}
\xyoption{all}

\usepackage{amsfonts,amsmath}

\begin{document}         

\newcommand{\C}{{\mathbb C}}
\newcommand{\e}{{\bf e}}
\newcommand{\F}{{\mathbb F}}
\newcommand{\G}{{\mathbb G}}
\newcommand{\Hom}{{\rm Hom}}            
\newcommand{\pee}{{\mathfrak p}}
\newcommand{\oh}{{\bf o}}
\newcommand{\pt}{{\rm pt}}
\newcommand{\Q}{{\mathbb Q}}
\newcommand{\rings}{{\rm rings}}
\newcommand{\ess}{{\mathbb S}}
\newcommand{\spec}{{\rm spec \;}}
\newcommand{\U}{{\rm U}}
\newcommand{\Z}{{\mathbb Z}}

\begin{center}
{\bf Complex cobordism and algebraic topology}
\end{center} \bigskip 

\dots it was my lot to plant the harpoon of algebraic topology into the
body of the whale of algebraic geometry \dots

S. Lefschetz, A page of mathematical autobiography, BAMS 74 (1968) 854 - 879
                            
Vers les ann\'ees 1962 m'arriva de rencontrer Lefschetz aux Etats-Unis. 
Il me fit alors cette confidence: ``Que l\'a Topologie \'etait belle avant
1935! Apr\`es, elle est devenue beaucoup trop alg\'ebrique." 

R. Thom, Probl\`emes recontr\'es dans mon parcours math\'ematique: un bilan,  
Publ.\ Math.\ IHES \ {\bf 70} (1989) \bigskip

\begin{center}
{\bf \S $\;$ Introduction}
\end{center}

This note is an attempt to supplement recent \begin{footnote} {This was written in 2004, and has
been minimally revised.} \end{footnote} [A2,Su2] accounts of Ren\'e Thom's work on
cobordism theory with a description of that subject's later evolution, under the influence of
Novikov, Quillen, Sullivan, and others.

Thom thought of himself above all as a geometer, but his thinking was 
grounded in formidable algebraic insight. Much of mathematics might be 
described as a dialogue between geometry and algebra (cf. e.g. Descartes), 
and although I will focus here on developments in {\bf algebraic} topology, 
I will tell this story in a lingua franca based on algebraic {\bf geometry}, 
not so far from the dialect used by Thom himself in his work [TL] on 
singularities of differentiable maps. I will try to justify the utility 
of that language in the course of a historical sketch. 

What follows is at best a tissue of generalities; it is a `soft' account, 
intended to be as accessible as possible. It therefore contains nothing 
about explicit calculations, which already have an extensive and accessible 
literature, cf. eg. [R1,2; WSW].

\begin{center}
{\bf \S 1 Preliminaries}
\end{center}

One of the most intimidating aspects of classical cobordism theory is the
sheer size of the rings involved. Thom showed that the unoriented
cobordism ring is graded polynomial, with one generator of each degree
not of the form $2^k - 1, \; k \geq 1$, and he showed that over the
rationals, the oriented cobordism ring is generated by the complex
projective spaces. Wall and others completed the calculation of the latter
ring, which has a considerable amount of two-torsion; but interest
eventually centered on the {\bf complex} cobordism ring $\Omega^\U_*$,
which was shown by Milnor and Novikov to be torsion-free, with one
generator in each even degree. We are therefore faced from the outset with
rings which are not finitely generated; but $\Omega^\U_*$ has the
technical advantage of being in some sense locally finitely presented (it is
{\bf coherent} [S]) which makes its homological algebra reasonably
tractable. Thus, although complex cobordism is in some ways less
geometrically natural than oriented cobordism, its algebraic accessibility
makes it a fundamental technical tool.

It took a while for attention to shift from the cobordism rings to their
associated (co)homology theories; indeed the multiplicity of such objects
played an important role in crystallizing our present notion of a spectrum.
That what we now call Thom spectra define cohomology theories was
used by Atiyah [A1] as a device to study relations between cobordism rings
associated to various families of classical groups, but before the work of
Novikov it was the cobordism rings {\it per se}, not their associated
cohomology theories, which were regarded as the significant subjects for 
investigation. Modern thinking tries to understand a cohomology theory
and its operations together, but with the technology available in the 60's,
this was an intimidating task [La]. Just as all cats look alike in the
dark, over the rationals all cohomology theories are equivalent: tensored
with $\Q$, the ring of cobordism operations must act (in some sense)
freely and transitively on the coefficients of the theory. The algebras of
cobordism operations are therefore approximately as enormous as the
coefficient rings themselves, and had Novikov not shown the remarkable
efficiency of an Adams spectral sequence based on complex cobordism as
a tool for calculating stable homotopy groups of spheres (and in particular
for clarifying the role played there by of the $J$-homomorphism [No]), our
understanding of these operations would probably have stayed limited for
a considerable time.

\begin{center}
{\bf \S 2 The geometric picture}
\end{center}
               
It was Quillen who gave us an economical language for these topics; but
he was very much concerned with making his paper as elementary and
self-contained as possible, so that language was backgrounded and largely
implicit. His formalism was strongly influenced by Grothendieck's
contemporary thinking about motives, and one way to paraphrase his
approach is to describe complex cobordism as the universal multiplicative
cohomology theory carrying good Chern classes for complex vector
bundles [Q; see also LM]. The work of Borel and Hirzebruch on
characteristic classes and symmetric functions reduces this to a question
about Chern classes for {\bf line} bundles, and the key issue becomes the
existence of a canonical (Euler, or first Chern) class $\e$ in the universal
example
\[
\Omega_\U^*(\C P^\infty) \cong \Omega_\U^*[[\e]] \;;
\]
thus any multiplicative cohomology theory $E^*$ with a Chern class for
line bundles is the recipient of a unique natural multiplicative
transformation $\Omega^*_\U \to E^*$ of cohomology theories.     

{\bf 2.1} Now the classifying space $\C P^\infty$ for complex line bundles
is a commutative $H$-space, so $\Omega^*_\U(\C P^\infty)$ acquires the
structure of a completed Hopf $\Omega^*_\U$-algebra on a single formal
generator $\e$. Conveniently enough, such `formal group laws' had been
extensively studied by Dieudonn\'e, Lazard, and others, as a theory of
Abelian arithmetic Lie groups. It is trivial to show the existence of a
universal one-dimensional formal group law (take the ring generated by
the coefficients of the Hopf diagonal, modulo the ideal of relations
imposed by associativity and commutativity), but Lazard demonstrated the
highly {\bf non}-trivial fact that the resulting ring $L$ is (implicitly
graded) polynomial over $\Z$, with one generator in each even positive
degree. Mi$\check{\rm s}\check{\rm c}$enko had previously identified
the formal logarithm for the complex cobordism ring as 
\[
\log_\Omega(\e) = \sum_{n \geq 1} \C P^{n-1} \; \frac{\e^n}{n} \;,
\]
and Quillen went on to show that the homomorphism
\[
L^* \to \Omega^*_\U
\]
classifying the formal group law defined by $\C P^\infty$ is in fact an
isomorphism: the formal group law of complex cobordism is {\bf
universal}. If only for psychological reasons, this is quite a remarkable
result; it provides the (enormous, amorphous) complex cobordism ring
with something very much like a preferred system of coordinates. 

[In fact explicit coordinates can be defined, at least locally with respect to
a prime $p$, by writing
\[
[p]_\Omega(\e) = p\e +_\Omega \cdots +_\Omega v_k \e^{p^k}
+_\Omega \cdots \;;
\]
where $[p]_\Omega$ denotes the effect on the formal group of multiplication 
by $p$ in the $H$-space structure on $\C P^\infty \sim H(\Z,2) \;, \; 
+_\Omega$ denotes addition with respect to the formal group law, and $v_k \in
\Omega^\U_{2(p^k-1)} \otimes \Z_{(p)}$ is a certain recursively defined
polynomial generator. These (Araki-Milnor) generators are very useful for
computations, but I will not pursue this question further: the purpose of
this note is really to sketch a {\bf coordinate-free} description of the
cobordism ring.]

{\bf 2.2.1} It is here that the language of algebraic geometry becomes
relevant. In commutative algebra one associates to a (commutative) ring
$A$, the space $\spec A$ of its prime ideals: such ideals $\pee$ have
integral domain quotients $A/\pee$, which imbed in their quotient fields, so
we can equally well think of an element of $\spec A$ as the kernel of some
ring homomorphisms from $A$ to a field $k$. In classical algebraic
geometry one restricted attention to {\bf maximal} prime ideals, which
correspond to points closed in the natural (Zariski) topology on $\spec
A$; in particular if $k$ is algebraically closed, and $A = \Z[x_1,\dots,x_n]
$ is polynomial, then such closed points correspond to assignments of
values in $k$ for the coordinate functions $x_k$, thus recovering classical
$n$-dimensional affine space over $k$. The language of schemes [AM]
provides a dictionary equating commutative rings to a category of
topological spaces endowed with a structure sheaf whose stalks are {\bf
local} rings; modules over these rings can thus be described as sheaves of
modules over this ringed space. 

In particular, a (graded-commutative) multiplicative cohomology functor
\[
 X \mapsto E^*(X)
\]
whose coefficient ring $E^*(\pt) = E^*$ is concentrated in even degrees
can without loss of generality be interpreted as taking values in a category
of (even-odd graded) sheaves over the ringed space $\spec E$ obtained by
ignoring the grading (for the moment; we'll return to this issue below).

{\bf 2.2.2} A more modern point of view, however, sees the prime ideal
spectrum as a special instance of a more general functor of points
\[
\spec A(-) := \Hom_\rings(A,-) : ({\rm Commutative} \; \rings) \to {\rm
(Sets)} \;,                   
\]             
the idea being that the classical prime ideal spectrum yields inconveniently
few points. Yoneda's lemma guarantees that the enriched functor $A
\mapsto \spec A(-)$ is faithful, ensuring the existence of an abundance of
points.

This interpretation developed out of the study of moduli problems, and it
fits perfectly with Lazard's result: if $k$ is a general commutative algebra,
then the elements of
\[
\Lambda(k) := \spec L(k):= \Hom_\rings (L,k)
\]                  
are precisely the formal group laws $F \in k[[X,Y]]$ such that
\[
 F(X,F(Y,Z) = F(F(X,Y),Z), \;  F(X,Y) = F(Y,X) = X
+ Y + \cdots \;,
\]
with coefficients in $k$. This allows us to describe $\spec L(k)$ as the set
of `all' (one-dimensional) formal group laws (i.e. over {\bf any}
convenient commutative ring). Lazard's theorem then describes this set of
group laws as an infinite-dimensional affine space $\Lambda$: specifying 
the particular ring we're working over is irrelevant, because this
formulation of the theorem behaves naturally under `change of rings'. 

{\bf 2.3.1} Quillen's theorem thus tells us that we can think of complex
cobordism as a cohomology theory taking values in the category of
sheaves of modules over the ringed space of `all' one-dimensional formal
group laws; moreover, when the functor is applied to a finite complex, the
resulting sheaf will be coherent. But picturesque as this interpretation may
be, it is perhaps not particularly compelling. It becomes more interesting,
however, when we consider cohomology operations.  
     
In the language of \S 2.2.2 above, a group object in the category of
schemes is precisely a representable functor from commutative rings to
groups; the additive group
\[
\G_a(A) = A = \Hom_\rings (\Z[t],A)
\]
and the multiplicative group
\[
\G_m(A) = A^\times = \Hom_\rings(\Z[t,t^{-1}],A)
\]
are typical (and important) examples. The operations in complex
cobordism can be described most easily in terms of the functor
\[
A \mapsto \Gamma(A) := \{ t(T) = \sum_{k \geq 0} t_k T^{k+1} \; | \;
g_0 \in A^\times \}
\]
which sends a commutative ring to the set of invertible formal power
series with coefficients in $A$. This set is in fact a group, with
composition $t,t' \mapsto t \circ t' := t(t'(T))$ as group operation; the
functor is representable (it is the spectrum, in the functorial sense, of the
polynomial ring 
\[
S = \Z[t_0,t_0^{-1}][t_k \; | \; k \geq 1] = S_*[t_0,t_0^{-1}] \;,
\]
and the diagonal
\[
\Delta t (T) = (t \otimes 1)((1 \otimes t)(T)) \in (S \otimes S)[[T]]
\]
makes $S$ into a commutative (but not cocommutative) Hopf algebra.
Thus $\Gamma$ is a groupscheme; it is easily seen to be an extension 
\[
1 \to \Gamma_0 \to \Gamma \to \G_m \to 1
\]
of the multiplicative group by iterated copies of the additive group, split
by a canonical imbeddding of the multiplicative torus $\G_m \to
\Gamma$. It acts naturally on the space $\Lambda$ of formal group
{\bf laws:} if $F \in \spec L(A)$ is a group law, and $t \in \Gamma(A)$ is
an invertible power series, then 
\[
F^t(X,Y) := t^{-1}(F(t(X),t(Y))) \in \spec L(A)
\]
is another formal group law. This action, called `change of coordinates',
is represented by a ring homomorphism
\[
L \mapsto L \otimes S                   
\]
satisfying the axioms in [Adams]; the resulting equivalence classes, or
orbits of the action of  $\Gamma$ on $\Lambda$, are called 
(one-dimensional) `formal groups'. The multiplicative subgroup of
$\Gamma$ plays a special role; it is sometimes useful to distinguish
between isomorphisms of formal group laws (represented by elements of
$\Gamma$) and {\bf strict} isomorphisms (represented by elements of
$\Gamma_0$).        

{\bf 2.3.2} The usual analysis of the operations on the (co)homology
theory associated to a spectrum {\bf E} requires the existence of a
K\"unneth theorem strong enough to guarantee that $E^*({\bf E})$ and
$E_*({\bf E})$ are in some sense dual Hopf $E^*$-algebras. This holds
for the complex cobordism spectrum {\bf MU}; in particular,
\[
\Omega^\U_*({\bf MU}) \cong L_* \otimes S_* 
\]
as (bilateral) Hopf algebras, from which it follows [Adams] that 
\[
\Omega^\U_*(X) = \pi_*(X \wedge {\bf MU})
\]
has a natural Hopf coaction
\[
\Omega^\U_*(X) \mapsto \Omega^\U_*(X) \otimes S_* \;.
\]
The role played by the grading in this brief summary has been suppressed,
and needs to be clarified. The category of comodules over the Hopf
algebra $S_0 := \Z[t_0,t_0^{-1}]$ representing the multiplicative group is
equivalent to the usual category of $2\Z$-graded modules, via the functor
which assigns to the evenly graded module $M_*$ its associated ungraded
module $M = \oplus_k M_{2k}$, with coaction
\[
\sum m_{2k} \mapsto \sum m_{2k} t_0^k : M \to M \otimes S_0 \;.
\]
The algebra $S_0$ acts on a formal group
\[
F(X,Y) = \sum F_{i,j} X^i Y^j 
\]
sending the coefficient $F_{i,j}$ to $t_0^{i+j-1}F_{i,j}$, thus endowing
$L$ with its intrinsic (even) grading; similarly, $t_k \in S_* \mapsto t_0^k
t_k$. But the cobordism $\Omega^\U_*(X)$ of a space will in general not
be concentrated in even degree, \begin{footnote}{Thanks to Haynes Miller for explaining to me
that such $\pm$-graded sheaves can be accommodated by adjoining a square root of $t_0$ to
$S_0$.}\end{footnote} so it is natural to regard it as an even-odd
graded comodule over $L \otimes S_0$, and to endow the tensor category
of such comodules with the Koszul associativity constraint usual in
algebraic topology [K]. 

{\bf 2.4.1} This discussion can now be summarized by saying that the
bordism module of a space $X$ can be interpreted, without loss of
generality, as an element $\Omega^\U_\bullet(X)$ of the tensor category
of even-odd graded sheaves of modules over the {\bf stack}
$\Lambda/\Gamma$ of one-dimensional formal groups: in other words, it
is a $\Gamma$-equivariant sheaf on the moduli space $\Lambda$ of
formal group laws. 

We can of course do homological algebra this category. The bordism
modules $\Omega^\U_\bullet (S^n)$ of the spheres define analogues of the
Tate twists in algebraic geometry, and Novikov's Adams spectral
sequence for stable maps between spaces $X$ and $Y$ has the groups
\[
{\rm Ext}^i_{\Lambda/\Gamma} (\Omega^\U_\bullet (X),
\Omega^\U_\bullet(S^kY))
\]
of extensions in this tensor category as its $E_2$ term. In particular, the
sheaf cohomology
\[
H^i_\Gamma(\Lambda, \Omega^\U_\bullet(S^k)) \Rightarrow \pi^S_{i-
k}(\pt) \;.
\]
of the moduli stack is the $E_2$ term of Novikov's spectral sequence for the
stable homotopy groups of the spheres. 

2.4.2 The stable homotopy category is the original, motivating example of 
what is now called a tensor triangulated category, and the derived
category of sheaves of modules over the stack $\Lambda/\Gamma$ is
another. Complex cobordism defines a homological functor on the first
category, and it is tempting to suspect that the Adams-Novikov spectral
sequence is a consequence of the existence of a lift of this homology
theory to a functor between these triangulated categories. It is significant 
that {\bf this is not the case:} a homology theory which lifts to the derived
category, in the sense that it takes distinguished triangles in the stable
homotopy category to distinguished triangles in some derived category of
modules, is necessarily ordinary homology [[CF]; but see, however, [Fr] and
[Ne]].

The derived category of sheaves of modules over the stack of one-dimensional 
formal groups is a remarkably good model for the stable homotopy category, 
but their precise relation is deep and is very far from being well-understood.
A great deal of research in recent years has been concerned with this 
question; the next section will describe some of that work. 

\begin{center}
{\bf \S 3 The orbit stratification}
\end{center}

It is a consequence of Lie theory that over a field of characteristic zero, 
all (one-dimensional) formal groups are isomorphic to the additive group 
$F_{\G_a}(X,Y) = X + Y$; this translates as the assertion that over such a 
field, the stack $\Lambda/\Gamma$ is equivalent to the category with one 
object, and the group $\G_m$ as morphisms. The category of even-odd graded 
sheaves over this quotient is thus equivalent (by H. Miller's categorical 
version of Shapiro's lemma [M2]) to the category of $\Z$-graded rational 
vector spaces. This is another reflection of the fact that as far as pure 
homotopy theory is concerned, all cohomology theories in characteristic zero 
are equivalent.   

In positive characteristic things are much more interesting. If $F$ is a
formal group law over a field $k$ of characteristic $p > 0$, it is easy to
see (e.g. by symmetric functions) that the operation of multiplication by
$p$ in the group structure is represented by a power series of the form
\[
[p]_F(T) = g(T^q)
\]
for some power $q = p^n$ of $p$ and some series $g$ with $g'(0) \neq
0$ (unless $[p]_F = 0$, when we take $n = \infty$ and say that $F$ is of
additive type). The integer $n$ is an invariant of the group {\bf law},
called its height; for example, the multiplicative formal group
$F_{\G_m}(X,Y) = X + Y + XY$ has
\[
[p]_{\G_m}(T) = (1 + T)^p - 1 \equiv T^p \; {\rm mod} \; p \;,
\]
and thus has height one. All positive integers occur: for example
Honda's logarithm 
\[
\log_{F_n}(T) = \sum_{k \geq 0} p^{-k} T^{q^k}    
\]
defines a group law over $\Z_{p}$ whose reduction modulo $p$ has
$[p]_{F_n}(T) = T^q$. 
                    
Height is a complete invariant for (one-dimensional) formal groups over a
separably closed field; in other words, in characteristic $p$ the set of
geometric points of $\Lambda$ stratifies into orbits indexed by the set
${\mathbb N} \cup \infty$. The orbit of a grouplaw $F$ is a homogeneous space
$\Gamma/\ess_F$, and an equivariant sheaf of modules over such a
quotient is {\bf locally free}; indeed, it can be recovered from the action of
the isotropy group $\ess_F$ of the orbit on the fiber of the sheaf above the
point $F$. This suggests the hope of understanding sheaves over the
moduli stack in terms of representations of the isotropy groups of the
orbits. These isotropy groups are very beautiful; they were the focus of
studies in the early sixties by several arithmeticians, and in this section I
will present some of their results. 

{\bf 3.1.0} The Chern class for line bundles in classical cohomology is
additive, which makes $n = \infty$ a reasonable example to begin with; but,
as we will see, formal groups of infinite height are in many ways
anomalous. The foregrounding of classical cohomology in most
presentations of algebraic topology may therefore to have distorted our
expectations of the nature of the stable homotopy category. Solving the
Hopf invariant one problem using classical cohomology requires all the
machinery of the Adams spectral sequence; this is in striking contrast to the
elegance and simplicity of Adams and Atiyah's proof using $K$-theory.
Ordinary cohomology lies on a stratum of infinite codimension in the
moduli space of complex-oriented theories: it may be simple in some ways,
but in others it is extremely degenerate. 

{\bf 3.1.1} The isotropy group for the additive group law over a field $k$
of characteristic $p$ assigns to a commutative $k$-algebra $A$, the group
\[
\ess_{\G_a}(A) = \{ t \in \Gamma(A) \;|\: t(X + Y) = t(X) + t(Y) \} \;;
\]
but in characteristic $p$ the function defined by $t(T) = \sum t_k
T^{k+1}$ is linear iff the coefficients vanish unless $k=p^s$ is a power of
$p$. A typical element of this group thus has the form
\[
a(T) = \sum a_k T^{p^k} \;;
\]
the composition
\[
(a \circ b)(T) = \sum a_k (\sum b_j T^{p^j})^{p^k} = \sum_i (\sum_{i = j
+ k} a_k b_j^{p^k}) T^{p^i} \;;
\]
of two such elements is represented by Milnor's formula
\[
\Delta \xi_i = \sum_{i = j + k} \xi_k \otimes \xi_j^{p^k}
\]
for the diagonal of the Hopf algebra $\F_p[\xi_k \;|\; k \geq 1]$ dual to
Steenrod's reduced $p$th powers. 

When $p=2$ this exhausts the operations in ordinary cohomology, but
there are more operations when the prime is odd. We will recover the
remaining operations below, but even without them it is no exaggeration
to paraphrase this little calculation by saying that the automorphisms of
the fiber of the cobordism sheaf, above the stratum of formal group laws
of additive type, are essentially the operations of ordinary cohomology. 

{\bf 3.1.2} The isotropy group $\ess_{\G_m}$ for the multiplicative group
law is much simpler. Multiplication by an integer in the formal group
structure defines an embedding
\[
n \mapsto [n]_{\G_m}(T) = (1 + T)^n - 1 : \Z \to {\rm End}_k(\G_m) \;,
\]
which extends to an embedding 
\[
\alpha \mapsto \sum_{k \geq 1} \binom{\alpha}{k} T^k \in k[[T]]
\]
of the ring $\Z_p$ of $p$-adic integers in the endomorphisms. The units
of this ring exhaust the automorphisms of the theory; they correspond to
the Adams operations in $p$-adically completed $K$-theory. 

It is important, and typical, that this groupscheme is essentially constant;
when $F$ has finite height over a field, the Lie algebra of $\ess_F$ is
trivial. For if $t$ is an automorphism of a grouplaw with $[p]_F(T)
= T^q$ (which, when $k$ is separably closed,  we can assume without loss
of generality), then $t \circ [p]_F = [p]_F \circ t$ implies that if $t(T) =
t_0T + \epsilon T^i$ (mod higher order terms) is an infinitesimal
automorphism (with $\epsilon^2 = 0$), then 
\[
t_0 T^q + \epsilon T^{iq} \equiv t_0^q T^q + \epsilon^q T^{iq} \equiv
t_0^q T^q 
\]
mod higher order terms; i.e. $\epsilon= 0$. In consequence, the
groupschemes $\ess_F$ are \'etale when the height of $F$ is finite; in
particular, if $F = F_n$ is the group law defined by the Honda
logarithm, its isotropy group $\ess_{F_n} = \ess_n$ can be
identified with the group of units of the noncommutative integral domain 
\[
\oh_n := W(\F_q) \langle F \rangle / (F^n - p) \;;
\]
here $W(k)$ denotes the ring of Witt vectors of $k$ (which can be
identified when $k =\F_q$ with the ring obtained by adjoining a
primitive $(q-1)$th root of unity to $\Z_p$), and the Frobenius element
$F$ (representing the endomorphism $F(T) = T^p$ of $F_n$) satisfies the
identity
\[
a^\sigma F = F a              
\]
for $a \in W(\F_q)$ (with $\sigma \in {\rm Gal}(\F_q/\F_p)$ the usual
$p$th power automorphism of the field with $q$ elements).

The isotropy groups of group laws of finite height thus become constant,
once the field of definition is large enough. Reduction modulo $p$
defines a ring homomorphism
\[
W(\F_q) \langle F \rangle /(F^n - p) \to \F_q \langle F \rangle /(F^n) 
\]
and thus a homomorphism on groups of units. The target is a truncation of the
group
\[
(\F_q \langle \langle F \rangle \rangle)^\times = \ess_{\G_a}(\F_q) 
\]
defined by the reduced power Hopf algebra; this leads to a consistency 
condition for the fibers of an equivariant sheaf over the moduli
stack. A finite complex defines a system of mod $p$ representations of the
groups $\oh_n^\times$ of strict units in a family of $p$-adic division
algebras, and these representations converge in a certain sense to a
representation of the Steenrod algebra. In intuitive terms: a finite complex
has limited ramification. 

{\bf 3.2.0} Quillen's theorem thus allows us to associate to a space, a
family of invariants associated to the stratification by height of the moduli
stack of formal groups. This stratification is in some ways quite simple: the
$n$th layer $\Lambda_n$ is the locus where the polynomial generators
$p = v_0, v_1, \cdots,v_{n-1}$ vanish, and $v_n$ is invertible. There is a big open
stratum, defined by formal groups of height one, and the fiber of the
cobordism sheaf above this generic orbit is (the mod $p$ reduction of)
classical $K$-theory; while at the bottom of the stratification we have a
closed orbit, carrying a sheaf with properties reminiscent of ordinary (mod
$p$) cohomology. 

At this point we are in the position of topologists in the days of Emmy
Noether, faced with the analog of torsion coefficients; we would like a
better understanding of the significance of these invariants. Ideas of
Sullivan permit us to interpret them in terms of closely related extraordinary
cohomology theories.

{\bf 3.2.1} To return for a moment to the wishful thinking in \S 2.4.2: if we
{\bf did} have a lifting
\[
({\rm Stable \; Homotopy}) \to D(\Lambda/\Gamma-{\rm Modules})
\]
of complex cobordism to a functor between tensor-triangulated
categories, then composition with the left-derived functor $Li^*_n$
defined by the inclusion morphism
\[
i_n : \Lambda_n/\Gamma \to \Lambda/\Gamma
\]
would define a functor taking values in equivariant sheaves over a
homogeneous space (and thus, after specializing at a typical fiber $F$, a
cohomology theory taking values in modular representations of $\ess_F$).
Such a derived functor would of course be more complicated than simple
restriction; in particular, it would imply the existence of a kind of 
universal coefficient spectral sequence
\[
E_2^{*,*} = {\rm Tor}^{\Omega_\U}(\Omega^\U_*(X),v_n^{-
1}\Omega^\U_*/I_n) \Rightarrow Li^*_n \Omega^\U_*(X)  \;,
\]
$I_n$ being the ideal $(p,v_1,\cdots,v_{n-1})$ defining the locus
$\Lambda_n$.  

It is quite remarkable that Sullivan's theory of cobordism with singularities
provides us with such a cohomology theory, taking values in equivariant
sheaves over $\Lambda_n$, calculable by just such a spectral sequence.
We cannot construct such a theory by applying homological algebra
naively to the coefficients, but we can mimic a standard construction in the
homotopy category. The sequence $v_0,\cdots,v_{n-1}$ of elements is
regular in the sense of commutative algebra, so the tensor product of the
complexes
\[
\xymatrix
{0 \ar[r]& \Omega^\U_* \ar[r]^{v_k}& \Omega^\U_* \ar[r]& 0}
\]
(for $0 \leq k \leq n - 1$) is a (Koszul) resolution of $\Omega^\U_*/I_n$.
In the homotopy category we can emulate this by taking the cofiber of a
suitable cubical diagram of spectra, constructed as a product of
elementary cofibrations realizing the elementary complexes above [M1].
The fiber of the resulting sheaf above the generic orbit is the mod $p$
reduction of complex $K$-theory, while the fiber above (a typical point
of) the closed orbit is ordinary (mod $p$) cohomology. More generally,
the fiber above a point of $\Lambda_n$ defines (the homological version of)
the `extraordinary $K$-theory' $K(n)^\bullet(X)$. 

Nowadays we have very elegant accounts [EKMM,HSS] of the stable category
which make these constructions relatively straightforward. I believe it 
is fair historically to say that the need to put these, and related, examples 
on a solid foundation provided much of the impetus which led to our current
highly polished understanding of the category of spectra. 

{\bf 3.2.2} It is not immediate from the description above that the cohomology
theories $K(n)$ are naturally $\ess_n$-representations; to prove this
requires us to exhibit a Hopf algebra structure on the operations, and for
that we need a rudimentary multiplicative structure [Ba,Ro]; there is a 
detailed account in [Wu]. The representing objects are definitely not 
$E_\infty$ ringspectra; very few mod $p$ theories are. It is however possible 
to identify their operation algebras completely (Boardman has completed the 
$p=2$ case only recently), and there is more to the story than the
isotropy groups. 

Sullivan's constructions use stratified sets with carefully controlled
singularities modelled on cones over a list of special example manifolds
[Su1]; in our case these models represent the classes $v_k$. Each such
model defines a generalized Bockstein operation on the resulting theory,
and in our case the operation defined by $v_k$ corresponds (in the sense
of \S 3.1.2) to Milnor's operation $Q_k$. The full Hopf algebra of operations
is thus an extension of some version of the group algebra of $\ess_n$ by an 
exterior algebra $E(Q_0,\cdots, Q_{n-1})$; to complete the description, we 
need to understand the action of the isotropy groups on the algebra of 
Bocksteins.

For this, a coordinate-free viewpoint is again useful. Lubin and Tate
constructed a cohomology theory $H^*(F,k)$ for formal groups over a
field $k$, and in particular showed that the second cohomology group
classified deformations. They calculate the rank of $H^2(F,k)$ to be the
height of $F$ minus one; an alternate formulation of their result identifies
this cohomology group (under its natural $\ess_F$-action) with the normal
bundle at $F$ of its orbit in $\Lambda$. We can thus think of the
Bocksteins as differentiation operators parametrized by deformations of the
group law. This accounts for all the Bocksteins except $Q_0$, which, in
this language, measures deformations in a direction away from $p$, toward
characteristic zero. 

\newpage

\begin{center}           
{\bf \S 4 Some more recent developments}
\end{center}

Stable homotopy theory allows us to understand spaces, which are
fundamental and perhaps irreducible objects of our imagination, in
intrinsically {\bf algebraic} terms. On the other hand, the deep formal 
properties of complex cobordism described above provide us with a {\bf
geometric} language for describing these objects (as sheaves of modules
over the structure space $\Lambda/\Gamma$). This intriguing double
reversal embeds algebraic topology into algebraic geometry in a way
which, it seems fair to say, was probably not anticipated by the creators of
the subject.

This concluding section collects some very brief remarks about more recent
work which takes for granted this circle of ideas, and pushes them forward.
It is idiosyncratic and incomplete, and in particular I will say little about 
elliptic cohomology, which is really the most fertile and striking new 
development in the subject. In the language presented above, the starting 
point of that theory is the idea that the formal completion of an elliptic 
curve at its origin defines a one-dimensional formal group. This defines a 
functor
\[
({\rm Elliptic \; curves}) \to ({\rm Formal \; groups}) \sim
\Lambda/\Gamma \;,
\]
or, more precisely, a morphism of {\bf stacks}. Roughly speaking, elliptic
cohomology is the pullback of complex cobordism, regarded as a sheaf
over the range of this morphism, to its domain; but this description
misrepresents the facts. For example, elliptic cohomology is {\bf not}
complex-oriented, and its construction (which is a tour-de-force [HGMR,Re])
requires a host of new ideas and techniques (see [Lu]) which I am not competent to
summarize. 

Elliptic cohomology is moreover a {\bf global} object: taking it apart one
prime at a time can be technically useful, but loses a great deal of
perspective. Its relations with conformal field theory on one hand, and the
equivariant topology of free loopspaces on the other, provide it with the
kind of rich connections to other fields which signal a major conceptual
advance [H2,AHS]. It is too early to attempt a summary; this story is just 
beginning. 

{\bf 4.1} The theory developed above derives about equally from algebra
(the theory of formal groups) and geometry (in particular, cobordism with
singularities). This cross-fertilization brings into focus the geometrically
unexpected existence of a {\bf hierarchy} of theories of singularity,
parametrized by the systems $p,v_1, \cdots, v_{n-1}$. There are in some
sense too many possible theories of cobordism with singularities, and it is
remarkable that algebra singles out this particular sequence of elements as
especially interesting.

Sullivan's constructions work with manifolds, not cobordism classes, and
the classes $v_n$ are represented (modulo their predecessors) by the $p$-dric 
hypersurfaces
\[
\{ [z_0: \cdots :z_q] \in \C P^q \;|\; \sum z_k^p = 0 \} 
\]
(where $q = p^n$). This suggests that the geometric significance of such
singularities should be sought in algebraic geometry in characteristic $p$;
but the question seems to have attracted little attention.

These models are extremely symmetric: they possess stationary-point-free
actions of $(\Z/p\Z)^n$ [F]. There is reason to think that the height
filtration of $\Lambda/\Gamma$ and the Atiyah-Swan filtration of the
spectrum of the cohomology of a finite group (by rank of supporting
elementary abelian groups) are two extreme aspects of some unified
phenomenon, whose understanding awaits progress in the study of
equivariant cobordism [HKR,GS,TCB]

{\bf 4.2} The deeper (arithmetic) aspects of the theory of formal groups
have been backgrounded in this account, through its use of the language of
algebraic varieties over fields. It is now clear that the moduli spaces of 
{\bf liftings} of formal groups (from characteristic $p$ to characteristic 
zero) [LT] provide natural sites for the localizations of stable homotopy 
theory, and work of Bousfield [B,BN,HvP] provides us with a systematic way to 
restrict to subcategories of spectra supported on suitably closed substacks of
$\Lambda/\Gamma$. In particular, the notion that the classical Euclidean
primes have chromatic overtones [H1] has become standard imagery. This is 
reflected by the deep, purely homotopy-theoretic fact [HS] that a the ring 
of stable endomorphisms of a finite complex has a center of Krull dimension 
one, generated either by some power of a prime $p$, or by a self-map which 
can be interpreted as a power of some $v_n$. This is reminiscent of the 
existence of a nontrivial element of the center of a finite $p$-group: it 
provides a natural, `internal' way for decomposing spaces inductively.

Our understanding of the relations between neighboring chromatic
primes is very primitive; the subject is really at the forefront of current
research. Waldhausen's chromatic red-shift principle [AR] suggests that
free loopspaces raise chromatic level; besides its connection with 
elliptic cohomology, this seems to have arithmetic echoes [AMS,TT]. 
Hopkins has proposed a chromatic splitting conjecture [Hv], which seems 
to require some further tinkering. There is now a very elegant
analysis [GHMR] of the sphere at height two (and $p=3$) which may serve
as a model for future research.

Underlying these developments looms the enormous question of the
homotopy-theoretic underpinnings of many of the constructions of classical
algebra. Quillen taught us that the $K$-theory of a ring, which before him
was purely a matter of generators-and-relations algebra, was in fact part of
homotopy theory; and that the $K$-groups themselves were really just the
homotopy groups of a much richer invariant defined by the $K$-theory
spectrum. Now we know that the moduli stack of formal groups is the tip
of another such homotopy-theoretic iceberg, and the list of such objects
has grown, to include the Lubin-Tate moduli spaces, the elliptic moduli
stack, and (in unpublished work of Hopkins) moduli spaces of $K3$
surfaces and more general Calabi-Yau manifolds. Waldhausen has
embraced this brave new world enthusiastically, and his point of view has
revolutionized algebraic $K$-theory. One is reminded of Hilbert's remark
about Cantor, that no one will expel us from the garden that he has shown
to us. 

{\bf 4.3} There is good reason [M3] to suspect that certain Shimura varieties [C2
\S 6.4], which play a central role in the proof of the local Langlands
conjecture, are ripe for such homotopification (see [BL]!). These objects parametrize
Abelian varieties decorated in various ways, e.g. with level structures as
generalized by Drinfel'd, and they lead to constructions with enormous
symmetry groups: roughly the product ${\rm Gl_n(\Q_p)} \times
D_n^\times \times {\rm Gal}(\overline{\Q}_p/\Q_p)$ (with $D_n$ the
division algebra defined by the field of quotients of the ring $\oh_n$ of \S
3.1.2). Such objects define a kind of Morita transformation relating
automorphic representations of ${\rm Gl}_n(\Q_p)$ to representations of
the unit groups of division algebras and to representations of Galois
groups. 

Work on the Langlands program has focussed on applications of this
construction to representation theory, but it is quite remarkable that {\bf
precisely} the same objects arise in algebraic topology in the context
sketched here [GH,HG,RZ]; though now it is the cohomology of these
objects (in the sense of \S 2.4.1), rather than their representations, which
is of interest. The discovery of these objects was one of the major
accomplishments of mathematics in the 20th century; understanding their
significance is a comparable challenge for the coming one. \bigskip

\begin{center}                          
{\bf References}
\end{center}

[Ad] J.F. Adams, {\bf Stable homotopy and generalised homology}, Chicago 
(1974, 1995)

[An] M. Ando, Isogenies for formal group laws and power operations in the
cohomology theories $E_n$, Duke Math. J. 79 (1955) 423 - 485

[AHS] -------, M. Hopkins, N.P. Strickland, Elliptic spectra, the Witten 
genus, and the theorem of the cube, Invent. Math. 146 (2001) 595 - 687

[AMS] M. Ando, J. Morava, H. Sadofsky, Completions of $\Z/(p)$-Tate 
cohomology of periodic spectra, Geometry \& Topology 2 (1998) 
145 - 174

[A1] M.F. Atiyah, Bordism and cobordism, Proc. Cam. Phil. Soc. 57 (1961) 
200 - 208

[A2] --------, The impact of Thom's cobordism theory, BAMS 41 (2004) 337 - 340

[AM] -----, I. MacDonald, {\bf Introduction to commutative algebra}, 
Addison-Wesley (1969)

[AR] C. Ausoni, J. Rognes, Algebraic $K$-theory of topological $K$-theory,
Acta Math. 188 (2002) 1 - 39

[Ba] A. Baker, $A_\infty$ structures on some spectra related to Morava 
$K$-theories. Quart. J. Math. Oxford 42 (1991) 403--419.

[Bf] A.K. Bousfield, On the Boolean algebra of spectra, Comm. Math. Helv. 54 
(1979) 368 - 377

[BL] M. Behrens, T. Lawson, {\bf Topological automorphic forms}, available at {\tt
arXiv:math/0702719}

[BN] M. B\"okstedt, A. Neeman, Homotopy limits in triangulated categories,
Compositio Math. 86 (1993) 209 - 234

[BCF] R.O. Burdick, P.E. Conner, E.E. Floyd, Chain theories and their 
derived homology, Proc. AMS 19 (1968) 1115 - 1118

[C1] H. Carayol, Nonabelian Lubin-Tate theory, p. 15 - 39 in {\bf Automorphic 
forms, Shimura varieties, and $L$-functions II}, Perspect. Math. 11, 
Academic Press (1990)  

[C2] ---------, Preuve de la conjecture de Langlands local \dots travaux de 
Harris-Taylor et Henniart, expos\'e 857 in {\bf Seminaire Bourbaki 1998/99}, 
Asterisque 266 (2000)

[D] V.G. Drinfel'd, Coverings of $p$-adic symmetric domains, Func. Anal. \&
Applications 10 (1976) 29 - 40 

[EKMM] A. Elmendorf, I. Kriz, M. Mandell, J.P. May, {\bf Rings, modules, and 
algebras in stable homotopy theory}, AMS Surveys and Monographs 47 (1997)

[F] E.E. Floyd, Actions of $\Z/p\Z)^n$ without stationary points, Topology 10 
(1971) 327 - 336

[Fr] J. Franke, Uniqueness theorems for certain triangulated categories with
Adams spectral sequences, available at {\tt www.math.uga.edu/archive.html}

[GH] B.H. Gross, M. Hopkins, Equivariant vector bundles on the Lubin-Tate 
moduli space, p. 23 - 88 in {\bf Topology and representation theory}, 
Contemp. Math. 158 AMS (1994)

[GHMR] P. Goerss, H-W Henn, M. Mahowald, C. Rezk, A resolution of the 
$K(2)$-local sphere,  Ann. of Math. 162 (2005) 777--822, available at 
{\tt arXiv:0706.2175} 
                              
[GS] J.D.C. Greenlees, N. P. Strickland, Varieties and local cohomology for 
chromatic group rings,  Topology 38 (1999) 1093 - 1139

[H1] M. Hopkins, Global methods in homotopy theory, p. 73 - 96 in 
{\bf Homotopy theory (Durham)}, LMS Lecture Notes 117, Cambridge (1987)

[H2] -----, Algebraic topology and modular forms, Proc. ICM I (Beijing 2002) 
291 - 317

[HGMR] -----, P. Goerss, H.R. Miller, C. Rezk, lectures at the 2003 M\"unster 
conference 

[HG] -------, B.H. Gross, The rigid analytic period mapping, Lubin-Tate space,
and stable homotopy theory, BAMS 39 (1994) 76 - 86

[HKR] ------, N.J. Kuhn, D.C. Ravenel, Generalized group characters and 
complex oriented cohomology theories, JAMS 13 (2000) 553 - 594

[HSD] ------, J.H. Smith, E. Devinatz, Nilpotence and stable homotopy theory 
I, Ann. Math 128 (1988) 207 - 241

[HS] ------, ------, Nilpotence and stable homotopy theory II, Ann. Math 148 
(1998) 1 - 49

[Hv] M. Hovey, Bousfield localization functors and Hopkins' chromatic 
splitting conjecture, p. 225 - 250 in {\bf The $\check{\bf C}$ech centennial} 
Contemp. Math. 181, AMS (1995)

[HvP] ------, J. Palmieri, The structure of the Bousfield lattice, in {\bf Homotopy invariant
algebraic structures} 175--196, Contemp. Math. 239, AMS (1999) 

[HSS] ------, B. Shipley, J.H. Smith, Symmetric spectra, JAMS 13 
(2000) 149 - 208

[HSt] ------, N.P. Strickland, {\bf Morava $K$-theories and localization}, 
AMS Memoir 666 (1999) 

[K] C. Kassel, {\bf Quantum groups}, Springer Graduate Texts 155 (1995)

[La] P. Landweber, Cobordism operations and Hopf algebras, Trans. AMS 129 
(1967) 94 - 110

[Lz] M. Lazard, Sur les groupes de Lie formels \'a un param\`etre, Bull. Soc.
Math. France 83 (1955) 251 - 274

[LM] M. Levine, R Pandharipande, Algebraic cobordism revisited, available at 
{\tt arXiv:math/0605196}

[LT] J. Lubin, J. Tate, Formal moduli for one-parameter formal Lie groups, 
Bull. Soc. Math. France 94 (1966) 49 - 59

[Lu] J. Lurie, A survey of elliptic cohomology, preprint

[MRW] H.R. Miller, D.C. Ravenel, W.S. Wilson, Periodic phenomena in the 
Adams-Novikov spectral sequence, Ann. Math. 106 (1977) 443 - 447

[M1] J. Morava, A product for the odd-primary bordism of manifolds with 
singularities, Topology 18 (1979) 177 - 186

[M2] --------, Noetherian localizations of categories of cobordism comodules,
Ann. Math 121 (1985) 1 - 39

[M3] --------, Stable homotopy and local number theory, in {\bf Algebraic analysis, geometry,
and number theory} 291--305, JHU (1989) 

[Ne] A. Neeman, Stable homotopy as a triangulated functor, Invent. Math. 109 
(1992) 17 - 40

[No] S.P. Novikov, Methods of algebraic topology from the point of view of 
cobordism theory, Izv. Akad. Nauk SSSR 31 (1967) 855 - 951

[Q] D.G. Quillen, Elementary proofs of some results of cobordism theory using
Steenrod operations, Adv. in Math. 7 (1971) 29 - 56

[R1] D.C. Ravenel, {\bf Complex cobordism and stable homotopy groups of 
spheres}, Academic Press Pure \& Applied Studies 121 (1986)

[R2] -------, {\bf Nilpotence and periodicity in stable homotopy theory}, 
Annals of Math. studies 128, Princeton (1992)

[Re] C. Rezk, Notes on the Hopkins-Miller theorem, p. 313 - 366 in {\bf 
Homotopy theory via algebraic geometry and group representations} (Evanston 
1997), 313--366, Contemp. Math. 220, AMS

[RZ] M. Rapaport, Th. Zink, {\bf Period spaces for $p$-divisible groups}, 
Annals of Math Studies 141, Princeton (1996)

[Ro] A. Robinson, Obstruction theory and the strict associativity of Morava 
$K$-theories, p. 143 - 152 in {\bf Advances in homotopy theory} (Cortona, 
1988), 143--152, LMS Lecture Notes 139 (Cambridge)

[S] L. Smith, On the finite generation of $\Omega^\U_*(X)$, J. Math. Mech. 18 
(1968/9) 1017 - 1023

[Su1] D. Sullivan, Singularities in spaces, p. 196 - 206 in {\bf Proceedings 
of Liverpool singularities II}, Springer Lecture Notes 209 (1971) 

[Su2] ----------, Ren\'e Thom's work on geometric homology and bordism, BAMS
41 (2004) 341 -350

[TL] R. Thom, Singularities of differentiable mappings, notes by H.I. Levine, 
in {\bf Proceedings of Liverpool Singularities I}, Springer Lecture 
Notes 192 (1971)

[TCB] C.B. Thomas, {\bf Elliptic cohomology}, Kluwer/Plenum (1999)

[To] T. Torii, On degeneration of one-dimensional formal group laws and 
applications to stable homotopy theory, Amer. J. Math. 125 (2003) 1037 - 1077

[Wa] F. Waldhausen, Algebraic $K$-theory of spaces, localization, and the
chromatic filtration of stable homotopy theory, in {\bf Algebraic topology, 
Aarhus}, p. 173 - 195 in Springer Lecture Notes 1051 (1984)

[W] W.S. Wilson, W, {\bf Brown-Peterson homology: an introduction and 
sampler}, CBMS Regional Conference Series in Mathematics 48 (1982)

[Wu] U. W\"urgler, Morava $K$-theories: a survey, in {\bf Algebraic topology:
Pozn\'an}, Springer Lecture Notes 1474 (1991)

\end{document}